\documentclass[12pt]{amsart}
\usepackage{amsmath, amsthm, amscd, amsfonts}
\newtheorem{Lemma}{Lemma}
\newtheorem{Theorem}{Theorem}

 \newfont{\Bbbb}{msbm10 scaled\magstephalf}
 \def\be{\begin{equation}}
 \def\ee{\end{equation}}

\def\bn{\begin{equation}}
\def\en{\end{equation}}
\def\br{\begin{center}}
\def\er{\end{center}}
\def\by{\begin{array}}
\def\ey{\end{array}}

\def\begy{\begin{eqnarray}}
\def\endy{\end{eqnarray}}
\def\bey*{\begin{eqnarray*}}
\def\eny*{\end{eqnarray*}}
\def\ber{\begin{tabular}}
\def\enr{\end{tabular}}

\def\bt{\begin{flushright}}
\def\et{\end{flushright}}

 \def\bea{\begin{equation}\begin{array}{ll}}
 \def\eea{\end{array}\end{equation}}

\begin{document}
\title[ composition operators]{ Composition operators  in the Lipschitz Space of the Polydiscs}
\author[Z.S. Fang and Z.H.Zhou]{Zhongshan Fang \and Zehua Zhou$^*$ }
\address{\newline Department of Mathematics\newline
Tianjin Polytechnic University
\newline Tianjin 300160\newline P.R. China.}

\email{fangzhongshan@yahoo.com.cn}

\address{\newline Department of Mathematics\newline
Tianjin University
\newline Tianjin 300072\newline P.R. China.}
\email{zehuazhou2003@yahoo.com.cn}

\keywords{ Composition operator, Lipschitz space, Polydiscs, Several
complex variables}

\subjclass[2000]{Primary: 47B38; Secondary: 26A16, 32A16, 32A26,
32A30, 32A37, 32A38, 32H02, 47B33.}

\date{}
\thanks{\noindent $^*$Zehua Zhou, corresponding author. Supported in part by the National Natural Science Foundation of
China (Grand Nos. 10671141, 10371091).}

\begin{abstract}In 1987, Shapiro shew that composition operator
induced by symbol $\varphi$ is compact on the Lipschltz space if and
only if the infinity norm of $\varphi$ is less than $1$ by a
spectral-theoretic argument, where $\varphi$ is a holomorphic
self-map of the unit disk. In this paper, we shall generalize
Shapiro's result to the $n$-dimensional case.
\end{abstract}

\maketitle


\section{Introduction}

Let $U^n$ be the unit polydiscs of $n$-dimensional complex spaces
$C^n$ with boundary $\partial U^n$, the class of all holomorphic
functions on domain $U^n$ will be denoted by $H(U^n)$. Let
$\varphi(z)=(\varphi_1(z),\cdots,\varphi_n(z))$  be  a holomorphic
self-map of $U^n$, composition operator is defined by
$$
C_\varphi(f)(z)=f(\varphi(z))$$ for any $ f\in H(U^n)$ and $z\in
U^n$.

In the past few years, boundedness and compactness of composition
operators between several spaces of holomorphic functions have been
studied by many authors: by Jarchow and Ried \cite{JR} between
generalized Bloch-type spaces and Hardy spaces, between Bloch spaces
and Besov spaces and BMOA and VMOA in Tian's thesis \cite{tja}.

More recently,  there have been many papers focused on studying the
same problems for $n$-dimensional case :  by Zhou and
Shi\cite{zs1}\cite{zs2}\cite{zs3} on the Bloch space in polydisk or
classical symmetric domains, Gorkin and MacCluer \cite{GorM} between
hardy spaces in the unit ball.

For the Lipschitz case, the compactness of $C_{\varphi}$ is
characterized by "little-oh" version of Madigan's \cite{Mad} the
bounedness condition, the same results in polydisc were obtained by
Zhou \cite{z4} and by Zhou and Liu \cite {zl}. In all these works
the main goal is to relate function theoretic properties of $\phi$
to boundedness and compactness of $C_{\phi}$.

To our surprise, by a spectral-theoretic argument, Shapiro
\cite{Sh2} obtained the following fact: $C_{\varphi}$ is compact on
the Lipschltz space $L_{\alpha}(D)$ if and only if
$||\varphi||_{\infty}<1$. In this paper, we shall generalize
Shapiro's result to the unit polydisc.

\section{Notation and background}

Throughout the paper, $D$ is the unit disk in one dimensional
complex plane, and $|||z|||=\max\limits_{1\leq j \leq n}\{|z_j|\}$
stands for the sup norm on the unit polydisc. Define $Rf(z)=<\nabla
f(z),\bar{z}>$ where $z=(z_1,\cdots, z_n)\in U^n$, and $H(U^n,D)$
for the class of the holomorphic mappings from $U^n$ to $D$. For
$0<\alpha <1$, it is well known that the Lipschitz space
$L_{1-\alpha}(U^n)$ is equivalent to $\alpha-Bloch$ space, which is
defined to be the space of holomorphic functions $f\in U^n$ such
that
$$||f||_{1-\alpha}=\sup\limits_{z\in U^n}\sum_{j=1}^n(1-|z_j|^2)^{\alpha}|\frac{\partial f}{\partial z_j}(z)|<\infty.$$
Here, Lipschitz space $L_{1-\alpha}(U^n)$ is a Banach space with the
equivalent norm $$||f||=|f(0)|+||f||_{1-\alpha}.$$

The Kobayashi distance $k_{U^n}$ of $U^n$ is given by
$$k_{U^n}(z,w)=\frac{1}{2}\log
\frac{1+|||\phi_z(w)|||}{1-|||\phi_z(w)|||},$$ where $\phi_z: U^n
\rightarrow U^n$ is the automorphism of $U^n$ given by
$$\phi_z(w)=(\frac{w_1-z_1}{1-\overline{z_1}w_1},\cdots,\frac{w_n-z_n}{1-\overline{z_n}w_n})$$
Since the map $t\rightarrow \log\frac{1+t}{1-t}$ is strictly
increasing on $[0,1),$ it follows that
$$k_{U^n}(z,w)=\max\limits_{1\leq j \leq n}\{\frac{1}{2}\log \frac{1+|\frac{w_j-z_j}{1-\overline{z_j}w_j}|}
{1-|\frac{w_j-z_j}{1-\overline{z_j}w_j}|}\}=\max\limits_{1\leq j
\leq n}\{\rho(z_j,w_j)\},$$ where  $\rho$ is the Poincar\'{e}
distance on the unit disk $D \subset C$.

Following \cite{Abate}, the horosphere $E(x,R)$ of center
$x\in\partial U^n$ and radius $R$ and the Kor\'{a}nyi region
$H(x,M)$ of vertex $x$ and amplitude $M$ are defined by
$$E(x,R)=\{z\in U^n: \limsup\limits_{w \rightarrow x}[k_{U^n}(z,w)-k_{U^n}(0,w)]<\frac{1}{2}\log R\}$$
and
$$H(x,M)=\{z\in U^n: \limsup\limits_{w \rightarrow x}[k_{U^n}(z,w)-k_{U^n}(0,w)]+k_{U^n}(0,z)<\log M\}.$$

We say that $f$ has $K-limit$ $L\in C$ at $x$ if $f(z)\rightarrow L$
as $z\rightarrow x$ inside any Kor\'{a}nyi region $H(x,M)$, we shall
write $\widetilde{K}-\lim\limits_{z\rightarrow x}f(z)=L$.

Let $f\in H(U^n,D)$ and $x\in \partial U^n$. If there is $\delta$
such that
$$\liminf\limits_{w\rightarrow x}\frac{1-|f(w)|}{1-|||w|||}=\delta
<\infty,$$ we call $f$ is $\delta-Julia$ at $x$. If there exists
$\tau \in  \partial U^n$ such that
$$f(E(x,R))\subseteq E(\tau,\delta R)$$ for all $R$,
we call this $\tau $ is the restricted $E$-limit of $f$ at $x$.

It should be noticed that $\delta>0$. In fact, $$\rho(0,f(w))\leq
\rho(0,f(0))+\rho(f(0),f(w))\leq \rho(0,f(0))+k_{U^n}(0,w);$$
therefore
$\frac{1-|f(w)|}{1-|||w|||}\geq\frac{1-|f(0)|}{2(1+|f(0)|)}>0$.

\section{Some Lemmas}

\begin{Lemma}  (Julia-Wolff-Carath\'{e}odory Theorem, Theorem 4.1 in \cite{Abate})
Let $f\in H(U^n,D)$ be $\delta-Julia$ at $x\in
\partial U^n$, and $\tau \in
\partial U $ be the restricted $E$-limit of $f$ at $x$, then
 $$\widetilde{K}-\lim_{z\rightarrow x} \frac{\partial f}{\partial x}(z)=\delta \tau.$$
\end{Lemma}

\begin{Lemma} (Theorem 1 in \cite{z4} or Corollary 4.1  in \cite{zl}) Composition operator $C_{\varphi}$
is bounded on the Lipschitz space $L_{1-\alpha}(U^n)$ if and only if
there is a constant $M>0$ such that
$$\sum\limits^n_{k,l=1}\left|\displaystyle\frac{\partial \phi_{l}}
{\partial z_k}(z)\right|
\left(\displaystyle\frac{1-|z_k|^2}{1-|\phi_l(z)|^2}\right)^{\alpha}\leq
M $$ for $z\in U^n.$\end{Lemma}

\begin{Lemma} (Theorem 2 in \cite{z4} or Corollary 4.2  in \cite{zl}) Composition operator
$C_{\varphi}$ is compact on the Lipschitz space $L_{1-\alpha}(U^n)$ if and only if
 $$\lim\limits_{\delta \rightarrow 0}\sup\limits_{dist(\varphi(z),\partial U^n)<\delta} \sum_{k,l=1}^{n}|
 \frac{\partial \varphi_l}{\partial z_k}(z)|\frac{(1-|z_k|^2)^{\alpha}}{(1-|\varphi_l(z)|^2)^{\alpha}}=0.$$\end{Lemma}

\begin{Lemma}(Lemma 3.2 in \cite{Abate}) Let $f\in H(U^n,D)$ and $x \in \partial U^n$. Then
$$\liminf\limits_{w\rightarrow
x}\frac{1-|f(w)|}{1-|||w|||}=\liminf\limits_{t\rightarrow
1^{-}}\frac{1-|f(\varphi_x(t))|}{1-t},$$ where $\varphi_x(z)=zx$ for
any $z\in D$.\end{Lemma}

\section{Main theorem}

\begin{Theorem} Suppose $C_\varphi$ is bounded on
$L_{1-\alpha}(U^n)$, then for every $ 1\leq l \leq n$ and $\xi \in
\partial U^n$ with $|\varphi_l(\xi)|=1$, $\varphi_l$ is
$\delta-Julia $ at $\xi$.\end{Theorem}

\begin{proof} For every $ 1\leq l \leq n$ and $\xi \in
\partial U^n$ with $\varphi_l(\xi)|=\eta$ and $\eta=e^{\theta_0}$, we will
show that $\varphi_l$ is $\delta-Julia $ at $\xi$ according to the
following cases.

{\bf Case 1:} $\xi=(\xi_1,\xi'), \xi_1=e^{\theta_1}$ and
$|||\xi'|||<1$.

First we consider the special case for $\xi=e_1=(1,0,\cdots,0)$ and
$\eta=1.$

For $ r\in (1/2,1)$ , define $\sigma(r)=(r,0,\cdots,0)=re_1$ such
that $$\lim\limits_{r\rightarrow 1^{-}} \varphi_l(\sigma(r))=1.$$
Setting $g(r)=\varphi_l(re_1)$, then $g'(r)=\frac{\partial
\varphi_l}{\partial z_1}(re_1)$. It follows from Lemma 2 that the
boundedness of $C_\varphi$ implies that
$$h(r)=R\varphi_l(re_1)(\frac{1-r}{1-\varphi_l(re_1)})^\alpha=rg'(r)(\frac{1-r}{1-g(r)})^\alpha$$
is bounded.

Putting $u(r)=\frac{1-g(r)}{1-r}$, it is easy to see that
$g'(r)=-(1-r)u'(r)+u(r)$ and
$$h(r)=ru(r)^{-\alpha}[-(1-r)u'(r)+u(r)].$$
If we write $v(r)=u(r)^{1-\alpha}$, then
$$-\frac{1}{1-\alpha}(1-r)v'(r)+v(r)=\frac{h(r)}{r}$$
the general solution of this differential equation is
$$v(r)=-\frac{1-\alpha}{(1-r)^{1-\alpha}}\int_1^r \frac{h(s)}{s(1-s)^{\alpha}}ds+\frac{C}{(1-r)^{1-\alpha}}.$$
Since $h$ is bounded, the first term in the right above is a bounded
function of $r$, and moreover $v(r)$ is of the order
$o(\frac{1}{(1-r)^{1-\alpha}})$ as $r\rightarrow 1^{-}$, so we have
$C=0$. Hence $v$, and moreover $u$ is also bounded, according to
Lemma 4, for some $\delta$, $\varphi_l$ is $\delta-Julia$ at $e_1$.

Now we return to the proof in case 1. Considering the mapping
$\tilde{\varphi}_l: U^n \rightarrow U^n,$ where
$$\tilde{\varphi}_l(z_1,z')=e^{-i \theta_0}\cdot \varphi_l
(e^{i\theta_1}z_1,\phi_{\xi'}(z')) $$ for $z=(z_1,z')\in U^n$. It is
easy to check that $C_{\tilde{\varphi}_l}$ is bounded on
$L_{1-\alpha}(U^n)$ and $\tilde{\varphi}_l(e_1)=1$.

By the above argument, we get $\liminf\limits_{t\rightarrow
1^{-}}\frac{1-|\tilde{\varphi}_l(te_1)|}{1-t}=\delta<+\infty$, that
is\begin{eqnarray*} \liminf\limits_{t\rightarrow
1^{-}}\frac{1-|\varphi_l(t\xi_1,\xi'))|}{1-t}&=&\liminf\limits_{t\rightarrow
1^{-}}\lim\limits_{r\rightarrow
1^{-}}\frac{1-|\varphi_l(t\xi_1,r\xi'))|}{1-t}\\
&\geq& \liminf\limits_{t\rightarrow
1^{-}}\frac{1-|\varphi_l(t\xi_1,t\xi'))|}{1-t}.
\end{eqnarray*}
It follows from Lemma 4 that $$\liminf\limits_{w\rightarrow
\xi}\frac{1-|\varphi_l(\xi)|}{1-|||\xi|||}=\delta <+\infty.$$

{\bf Case2:} $\xi=(\xi_1,\xi_2,\xi'), \xi_1=e^{\theta_1},
\xi_2=e^{\theta_2}$ and $|||\xi'|||<1$.

Now assume $\varphi_l(1,1,0,\cdots,0)=1$, and set
$g(r)=\varphi_l(r,r,0,\cdots,0)$ for $r\in (1/2,1)$. then
$g'(r)=\frac{\partial \varphi_l}{\partial
z_1}(r,r,0,\cdots,0)+\frac{\partial \varphi_l}{\partial
z_2}(r,r,0,\cdots,0)$, and so $R\varphi_l(r,r,0,\cdots,0)=rg'(r)$,
we can deal with it as in the case 1, and we can get $u$ is bounded,
furthermore $$\liminf\limits_{w\rightarrow
\xi}\frac{1-|\varphi_l(\xi)|}{1-|||\xi|||}=\delta <+\infty.$$

{\bf Case 3:} For the case $\varphi_l(\xi)=1$ with
$\xi=\sum\limits^n_{k=1} \beta_k e_k$, where $\beta_k=0$ or $1$, and
$e_k=(0,,0,\cdots,1,0,\cdots,0)$ with the $k-th$ component is $1$,
otherwise $0$; and even more general case, in a similar argument
with the cases 1 and 2, we can also show
$$\liminf\limits_{w\rightarrow
\xi}\frac{1-|\varphi_l(\xi)|}{1-|||\xi|||}=\delta <+\infty.$$ This
completes the proof of this theorem.
\end{proof}

\begin{Theorem} $C_{\varphi}$ is compact on
$L_{1-\alpha}(U^n)$ if and only if  $\varphi_j \in
L_{1-\alpha}(U^n)$ and $||\varphi_j||_{\infty}<1$ for each
$j=1,2,\cdots,n.$\end{Theorem}

\begin{proof} Sufficiency is obvious. Now we just turn to the necessity. Suppose
to the contrary
 that there exists $l$ ($1\leq l\leq n$) satisfying
$|\varphi_l(\xi)|=1$ for some $\xi \in \partial U^n$. It follows
from Theorem 1 that $\varphi_l$ is $\delta-Julia$ at $\xi$,
therefore by Lemma 1, we have $R\varphi_l(z)$ has $K-limit$ at
$\xi$. Hence
$$\sum_{k,l=1}^{n}|\frac{\partial \varphi_l}{\partial
z_k}(z)|\frac{(1-|z_k|^2)^{\alpha}}{(1-|\varphi_l(z)|^2)^{\alpha}}$$
\begin{eqnarray*}
&\geq& \sum_{k,l=1}^{n}|\frac{\partial \varphi_l}{\partial
z_k}(z)|\frac{(1-|||z|||^2)^{\alpha}}{(1-|\varphi_l(z)|^2)^{\alpha}}\\
&\geq&  \sum_{k,l=1}^{n}|z_k \cdot \frac{\partial
\varphi_l}{\partial
z_k}(z)|\frac{(1-|||z|||^2)^{\alpha}}{(1-|\varphi_l(z)|^2)^{\alpha}}\\
&\geq& C \sum\limits_{l=1}^n
|R\varphi_l(z)|\frac{(1-|||z|||)^{\alpha}}{(1-|\varphi_l(z)|)^{\alpha}}\\
&\geq& C\delta^{1-\alpha}
\end{eqnarray*}
as $z\rightarrow \xi$ inside any Kor\'{a}nyi region, where we can
take $C=\frac{1}{2^\alpha}$. It is a contradiction to the
compactness of $C_\varphi$ by Lemma 3. Now the proof of Theorem 2 is
completed.
\end{proof}

\end{document}